\documentclass[a4paper, 10pt]{article}

\usepackage[utf8]{inputenc}
\usepackage[T1]{fontenc}
\usepackage{lmodern}
\usepackage[UKenglish]{babel}
\usepackage{vmargin}
\usepackage{framed}
\usepackage{amsmath,amsthm,amssymb}
\usepackage{tikz}%
\usepackage[shortlabels]{enumitem}
\usepackage{csquotes} %
\usepackage[backend=biber,defernumbers=true,%
giveninits=true,isbn=false,doi=false,url=false,%
bibencoding=inputenc,style=alphabetic]{biblatex}
\newtheorem*{theorem*}{Theorem}
\newtheorem{theorem}{Theorem}

\newtheorem{corollary}{Corollary}[theorem]
\newtheorem{proposition}{Proposition}
\newtheorem*{conjecture*}{Conjecture}
\newtheorem*{definition*}{Definition}
\newtheorem*{question*}{Question}
\newtheorem*{notation*}{Notation}
\newtheorem*{observation*}{Observation}

\newcommand{\currentstatementname}{}
\newtheorem*{genericclaim}{\currentstatementname}
\newenvironment{citedstatement}[1]{\renewcommand{\currentstatementname}{#1}\begin{genericclaim}}{\end{genericclaim}}
\theoremstyle{definition}
\newtheorem*{remarks*}{Remarks}
\newtheorem*{remark*}{Remark}
\theoremstyle{remark}
\newtheorem*{claim*}{Claim}
\newtheorem{claim}{Claim}[proposition]
\newtheorem{notationinproof}[claim]{Notation}
\newtheorem*{remarkinproof*}{Remark}

\newenvironment{proofclaim}{\begin{proof}[Proof of Claim]}{\end{proof}}
\newcommand{\Aut}{\operatorname{Aut}}
\newcommand{\GL}{\operatorname{GL}}
\newcommand{\PGL}{\operatorname{PGL}}
\newcommand{\SO}{\operatorname{SO}}

\newcommand{\rk}{\operatorname{rk}}
\newcommand{\bC}{\mathbb{C}}
\newcommand{\bG}{\mathbb{G}}
\newcommand{\bK}{\mathbb{K}}
\newcommand{\bP}{\mathbb{P}}
\newcommand{\bR}{\mathbb{R}}
\newcommand{\bZ}{\mathbb{Z}}
\newcommand{\generated}[1]{\left\langle#1\right\rangle}
\newcommand{\defgenerated}[1]{\left\langle#1\right\rangle_{\mathrm{def}}}

\AtBeginBibliography{\small}

\title{The geometry of involutions in ranked groups with a \textsc{ti}-subgroup%
}
\author{Adrien Deloro\footnote{
Sorbonne Université, Institut de Mathématiques de Jussieu-Paris Rive Gauche, CNRS, Université Paris Diderot.
Campus Pierre et Marie Curie, case 247, 4 place Jussieu, 75252 Paris cedex 5, France}  ~and Joshua Wiscons\footnote{Department of Mathematics and Statistics. California State University, Sacramento. Sacramento, CA 95819, USA}}

\begin{document}

\maketitle

\setcounter{page}{1}
\renewcommand{\thepage}{\arabic{page}}

\begin{quote}
Solomon saith, \emph{There is no new thing upon the earth}. So that as Plato had an imagination, \emph{That all knowledge was but remembrance}; so Solomon giveth his sentence, \emph{That all novelty is but oblivion}. Whereby you may see that the river of Lethe runneth as well above ground as below.
\end{quote}
\begin{center}
\rule{.7\textwidth}{.5pt}
\end{center}

\abstract{We revisit the geometry of involutions in groups of finite Morley rank. The focus is on specific configurations where, as in $\PGL_2(\bK)$, the group has a subgroup whose conjugates generically cover the group and intersect trivially. Our main result is the subtle yet strong statement that in such configurations the conjugates of the subgroup may not cover all strongly real elements. As an application, we unify and generalise numerous results, both old and recent, which have exploited a similar method; though in fact we prove much more. We also conjecture that this path leads to a new identification theorem for $\PGL_2(\bK)$, possibly beyond the finite Morley rank context.}
\begin{center}
\rule{.7\textwidth}{.5pt}
\end{center}
\begin{center}
\S~\ref{S:introduction}. Introduction \quad --- \quad \S~\ref{S:geometric}. Theorem~\ref{t:geometric} \quad --- \quad \S~\ref{S:linear}. The B-Sides
\end{center}

\section{Introduction}\label{S:introduction}%
The present research belongs to the intersection of geometric algebra and model theory. Our work
follows two major themes in the theory of groups of Lie type, i.e.~the analysis of transformation groups of linear geometries:
\begin{itemize}
\item
the study of the distribution of involutions (elements of order 2), strongly real elements (products of two elements), and their centralisers, as made systematic by Brauer and his school and exemplified by the classification of the finite simple groups;
\item
the identification of \emph{inner} geometries, for instance incidence geometries which are already present ``inside'' the group (as opposed to representation theory where the geometry is some external, additional structure, or Chevalley theory where the Lie algebra and adjoint action are given).
\end{itemize}
Our context is model-theoretic, namely that of groups of finite Morley rank, generalising the familiar setting of algebraic groups over algebraically closed fields.

\subsection{A theorem and conjecture}

\paragraph{Two forms of one group and a conjecture.}
We investigate a family of configurations capturing two groups with much geometry in common: $\SO_3(\bR)$ and $\PGL_2(\bC)$.

Consider $\SO_3(\bR)$. The set $I$ of its involutions is exactly that of half-turns of the space; each $i \in I$ has an axis $L_i$ and a rotation plane $P_i = L_i^\perp$. Also note that distinct involutions $i$ and $j$ commute if and only if $L_i$ is contained in $P_j$ (and dually). Thus, introducing the relation $(i \preceq j) \iff (i j = ji \neq 1)$, one finds that the incidence structure $(I, I, \preceq)$ is naturally isomorphic to the projective plane $\bP^2(\bR) = (\{$vector lines$\}, \{$vector planes$\}, \leq)$. As such $\SO_3(\bR)$ as an abstract group \emph{defines} $\bP^2(\bR)$ as an incidence geometry, a point which begs for several remarks.
\begin{enumerate}
\item
The notion of \emph{definability} is \emph{first-order definability}, as made precise by model theory.
\item
As first realised by Hilbert, the incidence structure $\bP^2(\bR)$ itself defines (a field isomorphic to) $\bR$; hence the group $\SO_3(\bR)$ defines the field $\bR$, with all model-theoretic consequences, including that $\SO_3(\bR)$ does \emph{not} have finite Morley rank.
\item
Treating involutions both as lines and dual planes highlights that there is an underlying polarity of the incidence structure.
\item
Hilbert's coordinatisation theorem holds exactly of \emph{arguesian} projective planes, viz.~those satisfying the Desargues 10-point axiom. However, every projective plane which embeds into a projective $3$-dimensional space is always arguesian, so if possible, identifying a $3$-dimensional geometry is quite preferable to a $2$-dimensional one. This is the point of view we shall adopt in the sequel. Incidentally, $\SO_3(\bR)$ can indeed be endowed with the structure of a projective $3$-dimensional space (inducing the above projective $2$-dimensional structure on $I$), exactly as we define it in a more general setting below.
\end{enumerate}

We now give the group-theoretic description of the phenomena. Let $G = \SO_3(\bR)$ and $C \simeq \SO_2(\bR)$ be the subgroup of rotations with given axis, say $\bR e_1$. Letting $i$ be the half-turn around $\bR e_1$, we see that $C = C_G^\circ(i)$ (the connected component of the centraliser). Every half-turn $w$ with axis orthogonal to $\bR e_1$ also commutes with $i$. Moreover, the action of $w$ is a reflection of the plane $P_i = \mathrm{Span}(e_2, e_3)$, so conjugation by $w$ inverts $C$. Hence $C_G(i)$ is not connected, and one finds that $C_G(i) = C \rtimes \generated{w}$. 
Moreover, the set of conjugates of $C$ forms a partition of $G$. Indeed, since every non-trivial rotation has a unique axis and $G$ is transitive on $\bP^2(\bR)$, $C$ intersects its distinct conjugates trivially, and the conjugates cover $G$. In symbols, $G = \sqcup_{g \in G/N_G(C)} C^g$, and in words, \emph{$C$ is a \textsc{ti}-subgroup whose distinct conjugates cover $G$}. %

As a matter of fact and at the core of this paper, everything above \emph{almost} holds in $G = \PGL_2(\bC)$. Let $C$ be the group of (classes of) diagonal matrices, which contains a unique involution $i$. Note that $C$ is inverted by the anti-diagonal involution $w$, so as before, $C_G(i) = C \rtimes \generated{w}$. Here again, distinct conjugates of $C$ meet trivially, but now the resulting disjoint union only contains a generic subset of $G$, not all of it. That is,  $C$ is a \textsc{ti}-subgroup whose conjugates \emph{almost} cover $G$. And returning to the geometry of involutions, in $\PGL_2(\bC)$ the geometry is that of a \emph{generically} defined projective plane: generic (but not all) pairs of points are connected by a unique line, and dually. 

Of course, one should expect that $\SO_3(\bR)$ and $\PGL_2(\bC)$ have much in common: they are two forms of the simple algebraic group $A_1$ (over reasonable fields). We suspect that the above geometry of involutions, which exists (perhaps only in a generically defined form) for groups possessing an almost self-normalizing \textsc{ti}-subgroup, characterises decent forms of the root system $A_1$. The formal conjecture is a bit more modestly stated and only targets the ``characteristic not $2$'' case.

\begin{conjecture*}[``$A_1$-Conjecture'']
Let $G$ be a connected ranked group with involutions but with no infinite elementary abelian subgroup. Suppose that $G$ has a definable, connected subgroup $C < G$ which is \textsc{ti} and almost self-normalising.
Suppose that $C$ has even index in $N = N_G(C)$.
Then $G \simeq \PGL_2(\bK)$.
\end{conjecture*}

\paragraph{Ranked groups.}
It remains to explain our technical framework, and why $\SO_3(\bR)$ disappeared.
We will focus on \emph{groups of finite Morley rank}, also known as \emph{ranked groups}, which is a class extending groups of $\bK$-points of algebraic groups in algebraically closed fields and which was first discovered by model theorists. This gives us a well-studied context, but has the side effect of ruling out $\SO_3(\bR)$ (as mentioned above, a point that hinges on the fact that infinite fields of finite Morley rank are algebraically closed).

More precisely groups of finite Morley rank are groups in which the class of first-order definable sets bears a dimension function, called the \emph{rank}, which behaves like the Zariski dimension in algebraic geometry. 
Groups in this class have a lot in common with groups of $\bK$-points, in particular they possess:
\begin{itemize}
\item
an efficient notion of connectedness and of a connected component;
\item
an efficient notion of ``largeness'', called \emph{genericity}: subsets of maximal rank.
\end{itemize}
We wish to stress that the methodological framework forbids number-theoretic, topological, geometric (in the differential/functorial/infinitesimal sense of the term), and Lie-theoretic techniques.
It has however been conjectured by Cherlin and Zilber independently that infinite simple groups of finite Morley rank ought to be groups of $\bK$-points of simple algebraic groups in algebraically closed fields, a statement which is both deep and wide open.
Later it was realised that the analogy with the classification of the finite simple groups goes well beyond a mere similarity of statements to include many methods as well, suggesting that involutions, strongly real elements, and their centralisers are key notions in the domain. But quite interestingly, the research we describe here does not build on the bulk of this effort known as the ``Borovik programme for the classification of simple groups of finite Morley rank''. A rough summary of what we shall need from the theory is in \S~\ref{s:selfguide}; the monograph \cite{BNGroups} can serve both as a reference and an introduction to the topic.

For the present work, could one have elected another technical environment? This is unclear to us. The conjecture above does not aim at identifying an abstract group (which implies identifying a specific field as well), but a root system. Since $\SO_3(\bR)$ is the archetype of our configurations, one may wish to work in a context encompassing both algebraically closed and real closed fields, or more model-theoretically put, both ranked and $o$-minimal worlds. The latter also has useful notions of connectedness and largeness; moreover, the $o$-minimal analogue of the Cherlin-Zilber conjecture is already known. This suggests that one could have a future look at a common methodological denominator of both settings, for instance dimensional universes as developed in \cite{DWLinearisation}. But for the moment we shall content ourselves with finite Morley rank, at the cost of no longer admitting $\SO_3(\bR)$ in our scope.

\paragraph{A theorem.}

In this article, we prove Theorem~\ref{t:geometric} below, together with various corollaries given in \S~\ref{s:corollaries}. Theorem~\ref{t:geometric} is perhaps a bit subtle. We will add context, but the list of corollaries may be the best way to gauge its strength. 

The reader unfamiliar with ranked groups is encouraged to think of $G$ as a group of $\bK$-points over some algebraically closed field, which comes equipped with a dimension function (here called rank) but no a priori field nor topology. However, the ``characteristic'' is still quite relevant, which we approach  as follows.
A group $G$ is $U_2^\perp$ (this stands for: no $2$-unipotent subgroups) if it contains no infinite elementary abelian $2$-group; this may be safely read as working in characteristic not 2. The other case has been successfully dealt with in \cite{ABCSimple}.

We need two further definitions specific to the configuration we are studying. A subgroup $C < G$ is \textsc{ti} (this stands for: trivial intersections) if: $(\forall g)\ (C \neq C^g \to C\cap C^g = 1)$. Also, the subgroup $C$ is almost self-normalising if it has finite index in its normalizer, which, importantly, implies that the union of conjugates of $C$ has full rank in $G$.

\begin{theorem}[``The Geometric Theorem'']\label{t:geometric}
Let $G$ be a connected, $U_2^\perp$, ranked group with involutions. Suppose that $G$ has a definable, connected subgroup $C < G$ which is \textsc{ti} and almost self-normalising.

Then $\bigcup_{g \in G} C^g$ does not contain all strongly real elements of $G$.
\end{theorem}

Theorem~\ref{t:geometric} is a revisitation of the classical analysis by Nesin of specific $\SO_3(\bR)$-ish configurations~\cite{NNonsolvable} in a ranked context. It is however new and quite more general; we believe it clarifies the picture greatly as it focuses on involutions, not on subgroups.
Its name refers to the fact that such a configuration naturally defines a $G$-invariant point-line-plane incidence structure that is, in fact, a three-dimensional projective space, just like in $\SO_3(\bR)$.

Analogues of Theorem~\ref{t:geometric} could be, and may  well have been, considered in various familiar settings. For example, one may wonder if the theorem holds in the Lie context with the last sentence replaced by ``If $\bigcup_{g \in G} C^g$ contain all strongly real elements of $G$, then $G \simeq \SO_3(\bR)$'', but we are unaware of such a result. In the finite case, one could ask the following---notice the analogy with the classical Brauwer-Fowler ``$|G|^{1/3}$'' lemma, the proof of which deals with the worst case scenario where centralisers of involutions are disjoint.

\begin{question*}
Let $G$ be a finite group (simple of Lie type if necessary) with a \textsc{ti}-subgroup $C$. Suppose that $\sqcup_{g \in G/N_G(C)} C^g$ has cardinal $\geq (1-\varepsilon)|G|$. Get an upper bound for the number of involutions/lower bound for the order of the centraliser of an involution.
\end{question*}

 However, the present article may be seen as praise of \emph{generic} methods for infinite configurations, with loss of effectiveness as a counterpart; and we do not (nor intend to) address this question here.

\paragraph{Final remarks.}
Theorem~\ref{t:geometric} is a promising first step towards the $A_1$-conjecture, resolution of which (or more precisely a technical variant of it) would eliminate one of three critical, minimal cases of the Cherlin-Zilber conjecture left open by \cite{DJInvolutive}, namely ``CiBo$_2$''. Also, it should not be missed that the  reconstruction of projective geometries from decent forms of $\PGL_2(\bK)$ is also the topic of \cite{WAbstract}; however, in our case, completing a \emph{partial} (more specifically, generic) projective geometry into a \emph{genuine} one is not expected to proceed as in Weisfeiler's work but rather using model-theoretic tools.

Finally, however tempting and faithful to the spirit of the Cherlin-Zilber Conjecture, it would be immodest to drop the assumption on $N$ in the $A_1$-conjecture: in the special ``strongly embedded'' case (see Proposition~\ref{p:mainalternative}), we have no contradiction in sight for now. The worst possibility is of course when $N = C$: then one has to tackle ``generically Frobenius'' ranked groups, of which not enough is known, and of which CiBo$_1$ \cite{DJInvolutive} could well be an extraordinarily pathological case. %
\subsection{Corollaries}\label{s:corollaries}

The proof of Theorem~\ref{t:geometric} is in \S~\ref{S:geometric}.
Here we list some of its consequences, which range from classical to new. Corollaries~\ref{c:genuinecovering} and~\ref{c:Nesin} will be quickly derived in \S~\ref{s:corollariesofA}; Theorem~\ref{t:hereditary} and its Corollaries~\ref{c:linear} and~\ref{c:goodtori} will be proved in~\S~\ref{S:linear}. Difficulties are better appreciated in light of the fact that a ranked analogue of the Feit-Thompson ``odd order'' theorem is nowhere in sight, and might well not exist.

\begin{corollary}\label{c:genuinecovering}
Let $G$ be a connected, $U_2^\perp$, ranked group. Suppose that $G$ has a definable, connected subgroup $C < G$ whose conjugates partition $G$. Then $G$ has no involutions.
\end{corollary}

Corollary~\ref{c:genuinecovering} is entirely new. It was our starting point; then we realised that it was enough to capture all generic and strongly real elements to force genuine covering. %
Its setting is \emph{not} the same as Jaligot's ``Full Frobenius groups'' \cite{JFull}, where malnormality (viz.~$(\forall g)\ (g \notin C \to C \cap C^g = 1)$) was assumed, and absence of involutions proved after Delahan and Nesin.
It is a trivial matter to show that if in Corollary~\ref{c:genuinecovering} one supposes in addition that $C$ is soluble or contains no unipotent torsion, then it is self-normalising, hence malnormal.
Such configurations have been revived by Wagner \cite{WBad} as a generalisation of the ill-named ``bad groups'', of which more will be said below.

Modulo their non-trivial but standard local analysis, one immediately retrieves Corollary~\ref{c:Nesin}, a celebrated classic. (Announced in rank~$3$ by Cherlin \cite[\S5.2, Theorem~1]{CGroups} with a flawed proof, the result has an interesting story. In an unpublished preprint \cite{BInvolyutsii} Borovik foresaw the importance of geometrisation; the geometric line was simultaneously followed by Borovik and Poizat in a Soviet-published, much delayed article \cite[Theorem~5.4]{BPSimple}, and by Poizat in his self-published, fast-track, book \cite[\S~3.8]{PGroupes}; Nesin's approach in rank $3$ \cite{NNonsolvable} was different and independent, and to be later generalised by Corredor.)

\begin{corollary}[{in rank $3$, Nesin \cite[Main Theorem~(2)]{NNonsolvable}; in general, Borovik-Poizat \cite[Theorem~5.4]{BPSimple} or Corredor \cite[Theorem~1(5)]{CBad}}]\label{c:Nesin}
Let $G$ be an infinite simple ranked group whose proper, definable, connected subgroups are all nilpotent. Then $G$ has no involutions.
\end{corollary}

Parenthetically let us observe that if, in the previous corollary, only the definable, connected, \emph{soluble} subgroups are supposed to be nilpotent, then this is open. Recall that the \emph{Borel subgroups} are the definable, connected, soluble subgroups that are maximal as such.

\begin{conjecture*}
Let $G$ be a simple ranked group whose Borel subgroups are all nilpotent. Then $G$ has no involutions.
\end{conjecture*}

Such a configuration is sometimes referred to as a ``bad group''; the intriguing terminology has varied over the years, so one could consider refraining from using it. Also notice that what really matters in these tight configurations is the impossibility to apply the Schur-Zilber method (for defining fields) inside Borel subgroups. One might therefore include such hypothetical objects in the wider class of ``asomic groups'', i.e.~groups not interpreting an infinite field---always bearing in mind that other methods for constructing a field exist.%

What prevents us at present from proving the last conjecture is essentially that, as opposed to the algebraic world, Borel subgroups need no longer be conjugate in definable subgroups of $G$. This contrasts with the following result derived from Theorem~\ref{t:geometric} with only marginal extra work and in which \emph{hereditarily conjugate} means that every definable connected subgroup of $G$ enjoys Borel conjugacy. (We did our best to avoid the assumption, and could not.)

\begin{theorem}[generalises {Borovik-Burdges \cite[Theorem~3.8]{BBDefinably}}]\label{t:hereditary}
Let $G$ be a connected, $U_2^\perp$, ranked group in which all Borel subgroups are nilpotent and hereditarily conjugate. Then the (possibly trivial) $2$-torsion of $G$ is central.
\end{theorem}

From Theorem~\ref{t:hereditary}, we easily obtain Borovik and Burdges' ``linear theorem'', based on the preliminary analysis by Poizat \cite{PQuelques} (and later Mustafin \cite{MStructure}). Corollary~\ref{c:linear} was a significant (though unpublished) result; indeed, even in this ``definably linear'' context, which one would think much easier than the abstract ranked setting, the Cherlin-Zilber conjecture is open and there is no known analogue of the Feit-Thompson theorem there either. Our proof is fundamentally new. 

\begin{corollary}[{Borovik-Burdges, \cite[Main Theorem]{BBDefinably}}]\label{c:linear}
Let $\bK$ be a ranked field of characteristic~$0$ and $G \leq \GL_n(\bK)$ be a simple, definable, \emph{non Zariski-closed} subgroup. Then $G$ has no involutions.
\end{corollary}

We also found something we were not looking for.

\begin{corollary}\label{c:goodtori}
Let $G$ be a connected ranked group in which all Borel subgroups are good tori. Then the (possibly trivial) $2$-torsion of $G$ is central.
\end{corollary}

\begin{remark*}
The authors learnt shortly before first submitting that Alt\i{}nel, Berkman, and Wagner had also obtained a very special case of Corollary~\ref{c:goodtori}, when the Prüfer $2$-rank is $1$ \cite[Theorem~8]{ABWSharpTwo}. Their methods are also geometric but seem to focus on the plane of involutions (and ultimately Bachmann's theorem \cite[\S8.3]{BNGroups}) instead of the full $3$-dimensional geometry that is present.
\end{remark*}

\subsection{A self-guide to the present paper}\label{s:selfguide}

We close this introduction by reminding our reader what they will need to deem familiar in order to follow the main lines of our arguments. More specific references will be provided at technical moments.
\begin{itemize}
\item
Basic group-theoretic definitions such as a strongly real element \cite[\S~10.1]{BNGroups} and a strongly embedded subgroup \cite[\S~10.5]{BNGroups}.
\item
Some rough notion of what a group of finite Morley rank is; the unfamiliar reader should systematically use the analogy with groups of $\bK$-points in algebraically closed fields, in which the Zariski dimension ``ranks'' the constructible class \emph{but there is no functorial nor topological information}.

Among our fundamentals are the additivity property of rank \cite[\S~4.1.2]{BNGroups}, enabling basic rank computations, and also the notion of connectedness \cite[\S~5.2]{BNGroups} and why generic (i.e.~full rank) subsets of a connected group must intersect. Something \emph{almost} holds if it does up to taking connected components.

At times we also need the notion of the \emph{definable envelope} of a subset $X \subseteq G$, which we denote by $\defgenerated{X}$, as opposed to \cite[\S~5.5]{BNGroups} where it is called \emph{definable closure} and denoted $d(X)$.
\item
Very interestingly, classical and fundamental results such as the Chevalley-Zilber generation lemma (``indecomposability theorem'' \cite[Theorem~5.26]{BNGroups}), the Schur-Zilber linearisation lemma (``field theorem'' \cite[Theorem~9.1]{BNGroups}), and the analysis of groups of small Morley rank \emph{play no role at all}.
We would be glad to take it as an indication that our result goes beyond the mere theory of groups of finite Morley rank. But on the other hand we shall invoke Macintyre's theorem \cite[Theorem~8.1]{BNGroups} that an infinite field with a rank function is algebraically closed: something not true in other model-theoretic settings.

What really matters is torsion; torsion-lifting \cite[Exercise~11 on p.~93]{BNGroups}, basic divisibility properties \cite[Exercises~9 and 13 on p.~78]{BNGroups}, but also the following.
\item
Recall that a $U_2^\perp$ group is one not containing an infinite elementary abelian $2$-group. In the ranked case, Sylow $2$-subgroups (which are as always conjugate \cite[Theorem~10.11]{BNGroups}) then become toral-by-finite \cite[Corollary~6.22]{BNGroups}, i.e.~finite extensions of $(\bZ/2^\infty\bZ)^d$, where $\bZ/2^\infty\bZ$ is the Prüfer quasi-cyclic $2$-group. In a connected, $U_2^\perp$, ranked group \emph{with} involutions, the integer $d$ known as the Prüfer $2$-rank is non-zero \cite{BBCInvolutions}.
\item
Torality principles from \cite{BCSemisimple} are key features: in particular, in a connected, $U_2^\perp$ ranked group, every $2$-element lies in a $2$-torus.
\item
It is essential to know what a good or decent torus \cite{CGood} is and, in particular, the generic behaviour of the Cartan subgroups $C_G(T)$ and their conjugacy, as well as their connectedness \cite[Theorem~1]{ABAnalogies}.
\end{itemize}
If such tools, and a few more to be cited during the technical arguments, were available in the $o$-minimal context, Proposition~\ref{p:mainalternative} might carry over to there.
The crux of the argument is however not from model theory and relies on two deep theorems.%
\begin{itemize}
\item
A \emph{projective plane} can be arguesian (i.e.~satisfy the Desargues property) or not; a projective plane is arguesian iff it can be coordinatised by a skew-field (Hilbert); a \emph{projective $3$-space} is always arguesian. See \cite{HFoundations} (more accurate references in due time).
\item
Borel's fixed point theorem \cite[Theorem~21.2]{HLinear}: provided $\bK$ is an algebraically closed field, the group of $\bK$-points of a soluble, connected algebraic group acting on a complete algebraic variety must have a fixed point.
\end{itemize}

We note that the remarkable work by Bachmann and his school on the group-theoretic axiomatisation of Euclidean geometry (viz.~incidence+metric) plays no part in our approach. We insist that their work, however beautiful, belongs to metric geometry: methodologically speaking, it can therefore be no part of a ranked study, where fields are algebraically closed.

\section{Theorem~\ref{t:geometric}}\label{S:geometric}

\begin{citedstatement}{Theorem~\ref{t:geometric}}
Let $G$ be a connected, $U_2^\perp$, ranked group with involutions. Suppose that $G$ has a definable, connected subgroup $C < G$ which is \textsc{ti} and almost self-normalising.

Then $\bigcup_{g \in G} C^g$ does not contain all strongly real elements.
\end{citedstatement}

The present section is devoted to proving Theorem~\ref{t:geometric}, with a word on Corollaries~\ref{c:genuinecovering} and~\ref{c:Nesin} in \S~\ref{s:corollariesofA}.
A brief outline of the proof itself is as follows, in the same notation as the statement.
\begin{itemize}
\item
The main alternative (Proposition~\ref{p:mainalternative} of \S~\ref{s:mainalternative}) is highly restrictive. If $C$ is \textsc{ti} and self-normalising, then either $N_G(C)$ is strongly embedded, or $G$ shares many features of groups of type $A_1$ over $\bC$ or $\bR$; which covers the structure of the Sylow $2$-subgroup, the generic distribution of involutions, and even the generic nature of $G$.
\item
The assumption that $\bigcup_G C^g$ contains all strongly real elements actually rules out strong embedding, so Proposition~\ref{p:mainalternative} has much to say.
But as opposed to $\PGL_2(\bK)$, the configuration is now unviable: this is proved in \S~\ref{s:inconsistency}, Proposition~\ref{p:inconsistency}. To find a contradiction we introduce a projective geometry; this is remarkably accelerated by going $3$-dimensional with a built-in polarity, a strategy which must be compared with Nesin's work.
\end{itemize}

\subsection{The main alternative}\label{s:mainalternative}

\begin{notation*}
Throughout this subsection and the next, $G$ will be a connected, $U_2^\perp$, ranked group with involutions, and $C < G$ a definable, connected, \textsc{ti} and almost self-normalising subgroup.

We let $N = N_G(C)$ and $I \subseteq G$ be the set of involutions; $I \cdot I$ is then the set of strongly real elements of $G$.
\end{notation*}

\begin{proposition}[A French Nocturne]\label{p:mainalternative}
Let $G$ be a connected, $U_2^\perp$, ranked group with involutions. Suppose that $G$ has a definable, connected subgroup $C < G$ which is \textsc{ti} and almost self-normalising.

Then the following are equivalent:
\begin{enumerate}[(i)]
\item\label{i:Nnotstronglyembedded}
$N = N_G(C)$ is \emph{not} strongly embedded;
\item\label{i:N/Ceven}
$N/C$ has even order;
\item\label{i:SdihedralandCiC}
the Sylow $2$-subgroup is as in $\PGL_2(\bC)$ \emph{and} for $i \in I$ (the set of involutions), $C_G^\circ(i)$ is a conjugate of $C$;
\item\label{i:Cinverted}
$N = C \rtimes \bZ/2\bZ$ (inversion action);
\item\label{i:Cstronglyreal}
$C$ consists of elements which are strongly real in $G$ (alt.: $C \subseteq I \cdot I$);
\item\label{i:genericisstronglyreal}
the generic element of $G$ is strongly real (alt.: $I\cdot I$ is generic in $G$);
\item\label{i:genericstronglyrealisCartan}
every strongly real element \emph{generic as such} is in $\bigcup_G C^g$ (alt.: $\rk\left(I\cdot I \setminus \bigcup_G C^g\right) < \rk\left(I\cdot I\right)$);
\item\label{i:Cnotconnected}
for $i \in I$, $C_G(i)$ is \emph{not} connected;
\item\label{i:genericdiameter2}
if $i, j$ is a \emph{generic} pair of involutions, there is $k \in I$ commuting with both.
\end{enumerate}
\end{proposition}

\begin{remarks*}\
\begin{itemize}
\item
The reader may, like us, muse upon the fact that the first four conditions and \ref{i:Cnotconnected} provide \emph{structural} descriptions; whereas, the other four yield \emph{generic} information. Having such an equivalence is strong praise for generic methods in the study of groups of finite Morley rank.
\item
All are satisfied in $\PGL_2(\bK)$ (where $C \simeq \bK^\times$ is the algebraic torus, say the diagonal subgroup). Algebraically, $\SO_3(\bR)$ has the properties (where $C \simeq \SO_2(\bR)$ is the stabiliser of any non-zero vector in the natural representation), but its dimension function does not agree with the Morley rank. Return to the introduction for more details.

The $A_1$-conjecture is precisely that Proposition~\ref{p:mainalternative} characterises $\PGL_2(\bK)$ among groups of finite Morley rank.
\item
In $\PGL_2(\bK)$ and $\SO_3(\bR)$ it also happens that $I \cdot I = G$, i.e.~\emph{every} element is strongly real (in $\PGL_2(\bK)$ every element is either unipotent or semi-simple, and strongly real in either case, while in $\SO_3(\bR)$ every element is a rotation inverted by some half-turn of the space, hence strongly real). This property does not seem to follow from Proposition~\ref{p:mainalternative}, although it would be a consequence of the $A_1$-conjecture. And the latter does not look noticeably easier under the extra assumption that $I \cdot I = G$.
\item
Knowing the mere structure of the Sylow $2$-subgroup in \ref{i:SdihedralandCiC} does not seem to be enough. But \cite{DJInvolutive} provides a thorough analysis of $C_G^\circ(i)$ in $N_\circ^\circ$-groups, and we now say a (very optional) word about the CiBoes.

In CiBo$_1$, the (to our current knowledge, possibly \textsc{ti}) centraliser $C_G(i) = C_G^\circ(i)$ \emph{is} strongly embedded.
CiBo$_3$ is irrelevant here: neither the centraliser of an involution nor that of a maximal $2$-torus is \textsc{ti}.
Under the assumption that its connected centralisers $C_G^\circ(i)$ are genuinely \textsc{ti}, CiBo$_2$ is a model of Proposition~\ref{p:mainalternative}, making the $A_1$-conjecture non-trivial.
\item
Proposition~\ref{p:mainalternative} easily implies that $I\cdot C$ is generic in $G$. The converse may fail, as shown by the possibly \textsc{ti} configuration of type CiBo$_1$. To our surprise, genericity of $I \cdot C$ does not seem to be a relevant property.
\item
Non-connectedness of $C_G(i)$ implies $C_G^\circ(i) = C$, but here again the converse may fail. And there are no claims whatsoever on whether $N = C$ or not in the strongly embedded case.

\end{itemize}
\end{remarks*}

\begin{proof}
By an obvious rank computation, $\Gamma = \bigcup_{g \in G} C^g$ is generic in $G$. It follows that $C$ contains the connected centraliser of each of its non-trivial elements: if $1 \neq x \in C$, then $C_G(x) \le N_G(C)$ since $C$ is \textsc{ti}, so $C_G^\circ(x) \leq N_G^\circ(C) = C$ by almost self-normalisation.

Now, $C$ also contains a maximal decent torus of $G$: letting $T \leq G$ be one such, we know that $\bigcup_{g \in G} C_G(T)^g$ is generic \cite[Lemma~4]{CGood}. Since so is $\Gamma$, up to conjugacy there is $1 \neq x \in C \cap C_G(T)$, whence $T \leq C_G^\circ(x) \leq C$. Hence by torality \cite[Theorem~3]{BCSemisimple}, every $2$-element is in a conjugate of $C$.

Since $C$ is \textsc{ti}, an element $x \neq 1$ belongs to at most one conjugate of $C$, which we then denote by $C_x$. As we just saw $C_x$ is defined (at least) for $x$ a generic element of $G$ or a $2$-element. We stress that $C_x$ does not stand for $C_G^\circ(x)$, but for the unique conjugate of $C$ containing $x$; in the strongly embedded case equality could fail and all we know is $C_G^\circ(x) \leq C_x$.

There will be two main blocks: \ref{i:Nnotstronglyembedded}--\ref{i:genericstronglyrealisCartan} and \ref{i:Cnotconnected}--\ref{i:genericdiameter2}.
Actually we treat \ref{i:Cnotconnected} somehow separately, since the current status of generic use of the Borovik-Cartan polar decomposition in odd type is not fully clarified.
This is why we shall prove: [\ref{i:Nnotstronglyembedded}--\ref{i:genericstronglyrealisCartan}]$\Leftrightarrow$\ref{i:genericdiameter2}$\Leftrightarrow$\ref{i:Cnotconnected}. Of course \ref{i:Cinverted}$\Rightarrow$\ref{i:Cnotconnected}, but we want to stress that equivalence of \ref{i:Nnotstronglyembedded}--\ref{i:genericstronglyrealisCartan} with \ref{i:genericdiameter2} does \emph{not} rely on unpublished material.
The proof begins.
\begin{description}[font=\normalfont\itshape]
\item[\ref{i:Nnotstronglyembedded}$\Rightarrow$\ref{i:N/Ceven}.]
Let $g \in G\setminus N$ such that $N \cap N^g$ contains an involution $k$. If $N/C$ had odd order, then we would obtain $k \in C\cap C^g$, whence $C = C^g$ and $g \in N_G(C) = N$, a contradiction.
\item[\ref{i:N/Ceven}$\Rightarrow$\ref{i:SdihedralandCiC}.]
Suppose that $N/C$ has even order. Lifting torsion, let $\alpha \in N \setminus C$ be a $2$-element with least possible order, so that $\alpha^2 \in C$. However $\alpha \in C_\alpha \neq C$, so $\alpha^2 \in C_\alpha \cap C = 1$ and $\alpha$ is an involution. Now $C_\alpha \neq C$ implies that $\alpha$ inverts $C$. If $i \in C$ is any involution then $C = C_G^\circ(i)$; moreover $i$ normalises $C_\alpha \neq C$ so $i$ inverts $C_\alpha$, which proves that $C$ contains a unique involution. Hence the Prüfer $2$-rank is $1$, and since $\alpha \notin C$ normalises $C$, the Sylow $2$-subgroup of $G$ is as in $\PGL_2(\bC)$ (folklore; \cite[Proposition~27]{DJSmall} for reference).
\item[\ref{i:SdihedralandCiC}$\Rightarrow$\ref{i:Cinverted}.]
Consider a Sylow $2$-subgroup, isomorphic to $\bZ/2^\infty\bZ \rtimes \bZ/2\bZ$, with central involution $i \in C$ and some non-central involution $k$. Every involution in $C$ is toral in $C$ and the Prüfer rank is $1$, so $i$ is the only involution in $C$. In particular $C \neq C_k$, so $k$ inverts $C$. 
Hence $N \leq C_G(i) \leq N(C_G^\circ(i)) = N$, so by Steinberg's torsion theorem \cite{DSteinberg}, $N/C$ has exponent $2$. By the structure of the Sylow $2$-subgroup, $N = C \rtimes \generated{k}$.
\item[\ref{i:Cinverted}$\Rightarrow$\ref{i:Cstronglyreal}.]
Clear.
\item[\ref{i:Cstronglyreal}$\Rightarrow$\ref{i:genericisstronglyreal}.]
As we know, $\Gamma = \bigcup_G C^g$ is generic in $G$.
\item[\ref{i:genericisstronglyreal}$\Rightarrow$\ref{i:genericstronglyrealisCartan}.]
A strongly real element generic as such is then just a generic element of $G$; and as we know, $\Gamma$ is generic in $G$.
\item[\ref{i:genericstronglyrealisCartan}$\Rightarrow$\ref{i:Nnotstronglyembedded}.]
Suppose \ref{i:genericstronglyrealisCartan} and yet that $N$ \emph{is} strongly embedded. Let $S = I \cdot I$ be the set of strongly real elements and $\Sigma = S \cap \Gamma$. The assumption is that $\rk (S\setminus \Sigma) < \rk S$, which implies $\rk \Sigma = \rk S$ (the degree of $S$ is a priori unknown; of course once $S$ is proved generic, this will be settled at once).

By strong embedding, all involutions in $N$ are actually in $C$: this is because $N$ conjugates its involutions \cite[Theorem~10.19]{BNGroups} and $C = N^\circ$ contains one, say $i$, with connected centraliser $D = C_G^\circ(i) \leq C$. Then the set of involutions of $C$ is $I_C = i^N$, with rank $\rk C - \rk D$. And the set $S_C = I_C \cdot I_C$ of strongly real elements of $C$ has rank at most $2 \rk C - 2 \rk D$.

We first estimate $\rk \Sigma$.
If $1 \neq x = rs \in C$ is strongly real in $G$ (i.e.~with $r, s \in I$), then $r, s \in N$ so actually $r, s \in C$: hence $x \in S_C$. This shows:
\[\Sigma = \bigcup_G (S \cap C^g) = \bigcup_G S_{C^g} = \bigcup_G S_C^g,\]
whence $\rk \Sigma \leq \rk G - \rk N + \rk S_C \leq \rk G + \rk C - 2 \rk D$.

Now consider the following restriction of the product map:
\[\mu: I_C \times (I \setminus N) \to S.\]
We contend that it has trivial fibres, and that the image avoids $\Sigma$. The former uses essentially the same argument as the latter: if $ir \in C^g$ in obvious notation, then $i, r \in N^g$ which forces $C = C_i = C^g = C_r$, a contradiction.
Therefore:
\[\rk I_C + \rk (I \setminus N) = \rk \mu\left(I_C \times (I \setminus N)\right) \leq \rk(S\setminus \Sigma) < \rk S.\]
(Notice that we used $\rk (S\setminus \Sigma) < \rk S$, not only $\rk \Sigma = \rk S$.)

Thanks to $\rk(I\setminus N) = \rk I$ and $G$-conjugacy of involutions (\cite[Theorem~10.19]{BNGroups} again) we derive the following contradiction:
\[
\rk S > \rk C - \rk D + \rk G - \rk D \geq \rk \Sigma.\]

This shows that \ref{i:Nnotstronglyembedded}--\ref{i:genericstronglyrealisCartan} are equivalent.
\item[\ref{i:Cnotconnected}$\Rightarrow$\ref{i:genericdiameter2}.]
Whether $N$ is strongly embedded or not, we know by \cite[Theorem~10.19]{BNGroups}, resp.~the implication \ref{i:Nnotstronglyembedded}$\Rightarrow$\ref{i:SdihedralandCiC} (and torality principles), that all involutions are conjugate, say $I = i^G$. Since fibres of the map $g\mapsto i^g$ are translates of $C_G(i)$, there is a correspondence between the generic element of $G$ and the generic involution (genericity over $i$). Moreover $I^2$ has degree $1$.

Now if for the generic pair of involutions $(r, s)$ the definable envelope $\defgenerated{rs}$ (viz.~the smallest definable subgroup containing $rs$) is $2$-divisible, then for $g \in G$ generic over $i$, so is $\defgenerated{ii^g}$. Thus, the generic version of the Borovik-Cartan polar decomposition applies, showing that $C_G(i)$ is connected: a contradiction. (The decomposition we refer to can be derived using the tools of \cite[\S~5]{BBCInvolutions}, but the anxious reader should note that \ref{i:Cnotconnected} is \emph{not} used further in the present article.)
\item[\ref{i:genericdiameter2}$\Rightarrow$\ref{i:Cnotconnected}.]
Suppose that $C_G(i)$ is connected. Let $(i, j)$ be a generic pair of involutions. By \ref{i:genericdiameter2} there is $k \in I$ commuting with both, so $i \in C_G(k) = C_G^\circ(k) \leq C_k$, forcing $C_i = C_k = C_j$ likewise, violating genericity.

So \ref{i:Cnotconnected} and \ref{i:genericdiameter2} are equivalent. We finally establish the connection with \ref{i:Nnotstronglyembedded}--\ref{i:genericstronglyrealisCartan}.
\item[\ref{i:genericdiameter2}$\Rightarrow$\ref{i:Nnotstronglyembedded}.]
Let $(i, j)$ be a generic pair of involutions. By assumption there is $k$ commuting with both; say $C_k = C_i^g$. Now if $N$ \emph{is} strongly embedded, then $i \in C_i \cap N(C_k) \leq N(C_i) \cap N(C_i)^g$ implies that $g \in N(C_i)$, whence $C_i = C_k = C_j$ likewise, against genericity of $(i, j)$.
\item[\ref{i:Nnotstronglyembedded}$\Rightarrow$\ref{i:genericdiameter2}.]
We make free use of all properties \ref{i:Nnotstronglyembedded}--\ref{i:genericstronglyrealisCartan}.
Let $i$ be the involution of $C$. Consider the restriction of the multiplication map:
\[\mu: (N \setminus C) \times (I \setminus N) \to G.\]
The fibres are trivial since if $(k, r) \neq (\ell, s)$, in obvious notation, have the same image, then $1 \neq \ell k \in C$ is inverted by $r$: forcing $r \in N$, a contradiction. Since $\rk(N\setminus C) = \rk C$ and $\rk (I \setminus N) = \rk G - \rk C$, we find that $\mu$ is generically onto $G$. Moreover a generic element $(k, r) \in (N \setminus C) \times (I \setminus N)$ is mapped to one of $G$, and in particular there is $g \in G$ such that $kr \in C^g$. Then $k, r$ centralise $i^g$: the pair $(k, r)$ has the required property, which carries through conjugation to any generic pair.
\qedhere
\end{description}
\end{proof}

\subsection{The inconsistent configuration}\label{s:inconsistency}

We now venture into the inconsistent: as opposed to Proposition~\ref{p:mainalternative} which covered some existing groups, the next deals with pathological configurations.

\begin{proposition}%
\label{p:inconsistency}
Let $G$ be a connected, $U_2^\perp$, ranked group with involutions. Suppose that $G$ has a definable, connected subgroup $C < G$ which is \textsc{ti} and almost self-normalising.
Then the following are equivalent:
\begin{enumerate}[(i), start=10]
\item\label{i:CcoversII}
every strongly real element is in $\bigcup_G C^g$;
\item\label{i:diameter2}
if $i, j$ are any involutions, there is $k \in I$ commuting with both;
\item\label{i:CcoversG}
$G = \bigcup_G C^g$;
\item
$G$ does not exist.
\end{enumerate}
\end{proposition}

\begin{remarks*}\
\begin{itemize}
\item
Here again, $\SO_3(\bR)$ has all properties.%
\item
Using Proposition~\ref{p:mainalternative} it is a straightforward consequence of~\ref{i:CcoversG} that $G = I \cdot I$; the latter is however not inconsistent as it holds in $\PGL_2(\bK)$.
\end{itemize}
\end{remarks*}

\begin{proof}\setcounter{claim}{0}
Each condition clearly implies at least one in Proposition~\ref{p:mainalternative}, which we use freely.

\begin{claim}
We have equivalence of \ref{i:CcoversII}, \ref{i:diameter2}, and \ref{i:CcoversG}.
\end{claim}
\begin{proofclaim}\
\begin{description}[font=\normalfont\itshape]
\item[\ref{i:CcoversII}$\Rightarrow$\ref{i:diameter2}.]
Let $i \neq j$ be involutions. Then we may assume that $ij \in C$: so $i$ and $j$ normalise $C$ and centralise its unique involution.

\item[\ref{i:diameter2}$\Rightarrow$\ref{i:CcoversG}.]
We cheat a bit as we first prove \ref{i:CcoversII}.

Let $x = ij$ be any strongly real element. By assumption, there is an involution $k$ commuting with $i$ and $j$, so $x \in C_G(k)$. As we know from Proposition~\ref{p:mainalternative}, $x$ is in $C_G^\circ(k)$ or an involution, so in any case,  $x\in \bigcup_G C^g$.

We now prove that $I \cdot I \cdot I = I \cdot I$. Let $i, j, k$ be any  three involutions. By assumption, there is an involution $r$ commuting with both $i$ and $j$ and an involution $s$ commuting with both $r$ and $k$. 
We may assume that $r \neq i, j$ as otherwise $ij$ is an involution or 1, which implies that $ij\cdot k\in I \cdot I \cup I \subset I\cdot I $ (by Proposition~\ref{p:mainalternative}). In particular, $ij\in C_G^\circ(r)$. %

Now, if $s =r$, then $k\in C(r)$, and there are two options: $k=r$ or $k$ inverts $C_G^\circ(r)$. In the former case, $jk$ is an involution or $1$, forcing $i\cdot jk\in I \cdot I$.  In the latter case,  $ijk\in I\cup\{1\}$ since $ij\in C_G^\circ(r)$. So we are left to consider when $s \neq r$.

Of course, $ijk = ijs \cdot sk$. On one hand,  as $s$ inverts $ij\in C_G^\circ(r)$, $ijs$ is an involution or $1$; on the other hand,  the same is true of $sk$.

Hence $\generated{I} = I \cdot I$ is a definable subgroup containing the generic element by Proposition~\ref{p:mainalternative}: so $I \cdot I = G \subseteq \bigcup_G C^g$.
\qedhere\end{description}
\end{proofclaim}

From now on we suppose \ref{i:CcoversII}--\ref{i:CcoversG}. First we strengthen \ref{i:diameter2} as follows.

\begin{claim}\label{cl:involutiongraph}
Let $i \neq j$ be two distinct involutions. Then there is a unique involution which commutes with and is distinct from both.
\end{claim}
\begin{proofclaim}
If $i$ and $j$ do not commute then an involution commuting with both (this exists by~\ref{i:diameter2}) must be distinct from both. If $i$ and $j$ do commute then $ij$ is again distinct from both. We must now prove uniqueness. So let $i, j, k, \ell$ be four pairwise distinct involutions such that $k$ and $\ell$ commute with $i$ and $j$; notice how perfectly symmetric the configuration is; we are after a contradiction. The relevant portion of the commuting-involution graph is as follows---the involution $r$, which is introduced below, may in fact be equal to one of the others.
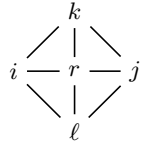
\begin{center}
\begin{tikzpicture}[line width = .6, scale = .8]
\node (i) at (-1,0) {$i$};\node (j) at (1,0) {$j$};
\node (l) at (0,-1) {$\ell$};\node (k) at (0,1) {$k$};
\node (r) at (0,0) {$r$};
\draw (i) -- (l) -- (j) -- (k) -- (i);
\draw (i) -- (r);\draw (j) -- (r);\draw (k) -- (r);\draw (l) -- (r);
\end{tikzpicture}
\end{center}

Let $r$ be the unique involution in $C_G^\circ(k \ell) = C_r$. We see that $i, j, k, \ell$ all normalise $C_r$, so as at most one of them can equal $r$,  we may suppose that $i, j, k \notin C_r$. Then by the structure of the Sylow $2$-subgroup, the fact that $\{1, i, k, r\}$ forms a four-group forces $i = kr$; likewise $j = kr = i$, a contradiction.
\end{proofclaim}

And now let us introduce an incidence geometry.

\begin{notationinproof}\label{n:Gamma}
Let $\Gamma$ be the incidence structure having:
\begin{itemize}
\item
for points the elements of $G$,
\item
for lines the translates of the various conjugates of $C$,
\item
for planes the translates of $I$,
\end{itemize}
and set-theoretic incidence relation (which is transitive). 
\end{notationinproof}

Notice that $G$ acts on the geometry $\Gamma$  both on the left and on the right---as opposed to if we had focused on the mere plane of involutions $I$, where only one copy of $G$ acts (by conjugation).

\begin{claim}
The action of $\bG = G \times G$ on $\Gamma$ is flag-transitive.
\end{claim}
\begin{proofclaim}
Indeed, $\bG$ is clearly transitive on planes, and as the stabiliser of the plane $I$ is the diagonal group $\bG^\Delta \simeq G$ (acting on $I$ as $G$ does by conjugation), this stabiliser is transitive on the lines of $I$, which are precisely those of the form $(N \setminus C)^g$. Consequently, the joint stabiliser of $I$ and $(N \setminus C)$ is $N$ in its conjugation action, which is easily seen to be transitive on the points of $(N \setminus C)$ by $2$-divisibility of $C$.
\end{proofclaim}

We also introduce a polarity, which will have the effect of considerably accelerating the proof.

\begin{notationinproof}
Let $\check{N} = N \setminus C$ and define an operation $\perp$ on $\Gamma$ as follows:
\begin{itemize}
\item
if $g$ is a point, let $g^\perp = g I$;
\item
if $L = gC$ is a line, let $(gC)^\perp = g \check{N}$;
\item
if $gI$ is a plane, let $(gI)^\perp = g$.
\end{itemize}
\end{notationinproof}

\begin{claim}
$\perp$ is a well-defined polarity of $\Gamma$, compatible with the action of $\bG$.
\end{claim}
\begin{proofclaim}
Clearly $\perp$ is well-defined for a point.
For a line, it amounts to having $C \check{N} = \check{N}$, an obvious claim.
Finally, for a plane: if $gI = hI$ then $x = h^{-1}g$ is strongly real, inverted by all involutions, so every involution normalises $C_x$, which we take to be $G$ if $x=1$. In fact, as we now see that $I\cdot I \subseteq N_G(C_x)$, we get $C_x = G$ and $x=1$. 

Of course $\perp$ is involutive on points and planes; on a line $L = gC$, fix $k$ inverting $C$ so that $\check{N} = kC$ and see that $L^{\perp\perp} = (g\cdot kC)^\perp = (gk C)^\perp = gk \cdot kC = gC = L$.

The polarity is clearly left $G$-covariant. It also is right-covariant: since $I$ is a normal set, only the case of a line could fail to be clear. But by left-covariance, this reduces to proving conjugacy-covariance, which is obvious.
It thus remains to check that the polarity preserves incidence. There are three cases to consider:
\begin{description}[font=\normalfont\itshape]
\item[point-line:]
by flag-transitivity and $\bG$-covariance, it suffices to prove it for $1 \in C$; now in usual notation, $C^\perp = N\setminus C \subseteq I = 1^\perp$;
\item[point-plane:]
here  we may work with $1 \in xI$, meaning that $x$ is an involution; hence $(xI)^\perp = x \in I = 1^\perp$;
\item[line-plane:]
translating again, suppose $C \subseteq xI$, forcing $x$ to be an involution inverting $C$; hence $(xI)^\perp = x \in \check{N} = C^\perp$.\qedhere
\end{description}
\end{proofclaim}

\begin{claim}
$\Gamma$ is a projective $3$-space.
\end{claim}

\begin{proofclaim}
The notion has a variety of equivalent definitions. We adopt the following: a projective $3$-space is a transitive point-line-plane incidence structure satisfying the following axioms as given by Hartshorne \cite{HFoundations}.
\begin{enumerate}[(S1)]
\item\label{threespace:i:twopoints} Every pair of distinct points lie on a unique line.
\item\label{threespace:i:threepoints} Three non-collinear points lie in a unique plane.
\item\label{threespace:i:lineplane} Every line and plane meet in at least one point.
\item\label{threespace:i:twoplanes} Two planes intersect in at least a line.
\item There exit four non-coplanar points, no three of which are collinear.
\item Every line has at least three points.
\end{enumerate}

First notice the following compatibility property: if $a \neq b$ are two distinct points of a plane $\Pi$ and $L$ is some line through $a$ and $b$, then $L \subseteq \Pi$. (This is a consequence of Hartshorne's axioms, but also something we need on our way towards them; notice that we make no claims on uniqueness of $L$ yet.)
By flag transitivity, it suffices to deal with $b = 1$ and $a \in C\setminus\{1\}$. Then a plane containing $1 \neq a$ must be of the form $k I$ where $k$ is an involution inverting $a$; notice that even if $a = i$ is the involution in $C$, then $k = i$ is forbidden. So $k$ inverts $C$, and $L = C \subseteq k I = \Pi$, as claimed.

We turn to proving the axioms.
The last two are trivial; Axiom~\ref{threespace:i:twopoints} is exactly $G = \bigsqcup_{G/N} C^g$, since for $x \neq y \in G$ only $x C_{x^{-1} y}$ will contain both. Axiom~\ref{threespace:i:twoplanes} is a weakening of \ref{threespace:i:twopoints}$^\perp$.
Axiom~\ref{threespace:i:threepoints} follows from another property.
\begin{quote}
(T) If a point does not lie on a line, then there is a unique plane containing both.
\end{quote}
(Existence in \ref{threespace:i:threepoints} follows from (T) modulo \ref{threespace:i:twopoints}; uniqueness also requires the ``compatibility'' above). Finally, Axiom~\ref{threespace:i:lineplane} is a weakening of the dual of (T), which we now prove in our setting.

By flag transitivity, we may suppose that $g \in G$ does not lie on $L = C$. Let $C_g = C_G^\circ(g) \neq C$ and $i_g \neq i$ be the involution in $C_g$. As we know from Claim~\ref{cl:involutiongraph}, there is an involution $j$ commuting to both $i$ and $i_g$ and distinct from both. Then $j$ normalises $C$ and $C_g$ but lies in neither, so it inverts both.  In particular, $g \in jI$ and $C \subseteq j I$, so the plane $jI$ meets the requirements.
We also contend that this plane is unique. If $xI$ contains $C$ (hence $1$), then $x$ is an involution, and if $xI$ also contains $g$, then it contains $i_g$ since this lies on the line through $1$ and $g$. Thus, $x$  centralises, $i$ and $i_g$, so uniqueness of the plane amounts to uniqueness of $j$ above, also given by Claim~\ref{cl:involutiongraph}.
\end{proofclaim}

\begin{remark*}
The Veblen-Young axioms for a projective space provide an alternative to the (disappointingly asymmetric) axioms of Hartshorne above. To use the former, one needs to show that if $a, b, c, d$ are distinct points for which the line through $a$ and $b$ meets the line through $c$ and $d$, then the line through $a$ and $c$ meets the line through $b$ and $d$.

Using the above ``compatibility'' property together with (T), this is rather trivial. Indeed, we may assume $a, b, c$ are noncollinear, so they lie in a plane, which we may take to be $I$. By assumption, the line through $c$ and $d$ has two points in $I$, so the entire line, hence $d$, is in $I$. Thus, what we are really showing is that two lines in $I$ must meet, which is just expressing that two distinct involutions commute with a third. 

But of course this alternative approach still uses our accurate knowledge of involutions.
\end{remark*}

\begin{claim}
Contradiction.
\end{claim}
\begin{proofclaim}
$G$ acts regularly by left-translation on the points of $\Gamma$. But $\Gamma$ is a $3$-dimensional projective space, hence coordinatisable as $\Gamma \simeq \bP^3(\bK)$ by the ``Desargues-Hilbert'' theorem (\cite[\S\S~24 sqq.]{HGrundlagen}, \cite[\S~2]{AGeometric}, or \cite[\S~7]{HFoundations} for the planar case; bear in mind that every $3$-dimensional projective space is arguesian, for instance \cite[Theorem~2.1]{HFoundations}).
Now $G$ acts regularly on this, implying $G \leq \Aut(\bP^3(\bK)) \simeq \PGL_{4}(\bK) \rtimes \Aut(\bK)$ \cite[Fact~3.5]{BNGroups}. Of course $\bK$ is algebraically closed by Macintyre's theorem on fields \cite[Theorem~8.1]{BNGroups} since Desargues field interpretation is definable (something obvious with Hilbert's proof  but harder to see with Artin's method---unless one understands Hrushovski's group configuration theorem); since there are no definable groups of field automorphisms in our ranked context \cite[Theorem~8.3]{BNGroups}, one finally finds $G \leq \PGL_4(\bK)$.

Consider the Zariski-closure $\overline{C}$ of $C$ inside $\PGL_4(\bK)$; being closed, connected, 
 and soluble (the latter can be seen as an instance of \cite[Corollary~5.38]{BNGroups}), it has a fixed point in its action on the complete variety $\bP^3(\bK)$ \cite[Theorem~21.2]{HLinear}. Hence so does $C \leq \overline{C}$: violating regularity.
\end{proofclaim}
Theorem~\ref{t:geometric} %
is now proved.
\end{proof}

We wish to make some more comments on the final contradiction.

\begin{remarks*}\
\begin{itemize}
\item
One does not need the full strength of Borel's theorem: since $C$ is abelian and the action of $\PGL_4(\bK)$ on $\bP^3(\bK)$ is a factor of the natural representation of $\GL_4(\bK)$, we may trigonalise simultaneously and find a fixed point for the preimage of $C$ in $\GL_4(\bK)$. So all we use geometrically seems to be the $2$-nilpotent version of the Lie-Kolchin theorem.
\item
It is well-known that Hilbert interpretation does not rely on Chevalley-Zilber indecomposable generation (parenthetically said, and contrary to widespread belief, Schur-Zilber field interpretation does not either: see \cite{DWLinearisation}), so definable coordinatisation would carry to other model-theoretic settings, say finite-dimensional theories.
\item
On the other hand, to find a fixed point one \emph{does} need algebraic closedness of $\bK$, which is typical of $\aleph_1$-categorical behaviour. And indeed, real closed fields do not satisfy Borel's theorem; and indeed, $\SO_3(\bR)$ does exist.

In the $o$-minimal case, one may argue in favor of using Bachmann's theorem in order to explicitly identify a real form of $\PGL_2$; but not in finite Morley rank in order to prove a contradiction.
\end{itemize}
\end{remarks*}

\subsection{Corollaries of Theorem~\ref{t:geometric}}\label{s:corollariesofA}

We shall be brief.

\begin{citedstatement}{Corollary~\ref{c:genuinecovering}}
Let $G$ be a connected, $U_2^\perp$, ranked group. Suppose that $G$ has a definable, connected subgroup $C < G$ whose conjugates partition $G$. Then $G$ has no involutions.
\end{citedstatement}
\begin{proof}
$C$ is \textsc{ti} and almost self-normalising by a clear rank computation; if $G$ contains involutions then Theorem~\ref{t:geometric} applies and $G = \bigcup_{g \in G} C^g$ does not contain all its strongly real elements, a contradiction.
\end{proof}

\begin{citedstatement}{Corollary~\ref{c:Nesin}}
Let $G$ be an infinite simple ranked group whose proper, definable, connected subgroups are all nilpotent. Then $G$ has no involutions.
\end{citedstatement}
\begin{proof}
We first invoke \cite{ABCSimple}: an infinite simple ranked group of \emph{infinite $2$-rank} is algebraic, hence has non-nilpotent Borel subgroups. As such, we may assume that $G$ is $U_2^\perp$. We are not using this major result out of mere laziness; we actually believe that \cite{ABCSimple} covers (and exhausts) a topic other than ours.

We then refer to the initial analysis in \cite[Theorem~13.3]{BNGroups}; there ads (i)--(iii) (or in the proof, Claims (a)--(e)) are straightforward and do not rely on much. In particular, we shall admit that all Borel subgroups (which here are the definable, connected, proper subgroups maximal as such) are conjugate and partition $G$.
Adding involutions to the picture brings us under the assumptions of Theorem~\ref{t:geometric}, and this contradicts the partitioning property.
\end{proof}

\section{The B-Sides}\label{S:linear}

Here is the other theorem we announced in the introduction. We remind the reader that \emph{hereditarily conjugate} means that every definable connected subgroup of $G$ enjoys Borel conjugacy, viz.~conjugacy of its maximal definable, connected, soluble subgroups.

\begin{citedstatement}{Theorem~\ref{t:hereditary}}
Let $G$ be a connected, $U_2^\perp$, ranked group in which all Borel subgroups are nilpotent and hereditarily conjugate. Then the (possibly trivial) $2$-torsion of $G$ is central.
\end{citedstatement}

We shall first prove Theorem~\ref{t:hereditary}, then its Corollaries~\ref{c:linear} and~\ref{c:goodtori}, which will be recalled hereafter.

\subsection{Proof of Theorem~\ref{t:hereditary}}

\setcounter{claim}{0}
\begin{proof}
Let $G$ be a counterexample of minimal rank. Notice that any definable, connected subgroup of $G$, or any quotient by a definable, normal, \emph{soluble} subgroup of $G$ still satisfies the assumptions (in the latter case, the rank need not be smaller).

\begin{claim}
We may assume that $G$ is centreless.
\end{claim}
\begin{proofclaim}
Suppose the theorem is proved for centreless groups.

If $Z(G)$ is finite, then $G/Z(G)$ satisfies all hypotheses \emph{and} is centreless \cite[Lemma~6.1]{BNGroups}; by our extra assumption, the $2$-torsion of $G/Z(G)$ is central, hence trivial as $Z(G/Z(G)) = 1$: this means that the $2$-torsion of $G$ was in $Z(G)$, as desired.

If $Z(G)$ is infinite, then $G/Z(G)$ satisfies all hypotheses but has smaller rank; by induction the $2$-torsion of $G/Z(G)$ is central, 
meaning that the $2$-torsion of $G$ is in $Z_2(G)$ (the second centre).
By \cite[Theorem~6.9]{BNGroups}, the (therefore unique) maximal $2$-torus $T$ of $G$ is central in $Z_2(G)$, hence characteristic in $G$, hence central in $G$ \cite[Theorem~6.16]{BNGroups}. So $T \leq Z(G)$ and $G/Z(G)$ has no non-trivial $2$-torus, which by \cite[Theorem~1]{BBCInvolutions} implies that it actually has no involutions. Here again all the $2$-torsion of $G$ was in $Z(G)$.
\end{proofclaim}

Suppose that $G$ contains involutions; let $T > 1$ be a maximal $2$-torus of $G$, and $C = C(T)$ be its centraliser, which we know by \cite[Theorem~1]{ABAnalogies} to be connected. Also let $N = N(T) = N(C(T))$ be its normaliser, a finite extension of $C$ by the rigidity of tori \cite[Theorem~6.16]{BNGroups}.

\begin{claim}
$C$ is a \textsc{ti}-subgroup of $G$.
\end{claim}
\begin{proofclaim}
Let $C_1, C_2$ be conjugates of $C$ meeting in $x \neq 1$. Then $H = C_G^\circ(x)$ is a definable, connected, proper subgroup of $G$; by induction, its $2$-torsion is central. Since it contains both $T_1$ and $T_2$ (the maximal $2$-torus of $C_1$, resp.~$C_2$), this forces $T_1 = T_2$ and $C_1 = C_2$.
\end{proofclaim}

\begin{claim}
If $x \in G$ is a strongly real element then $x^2 \in \bigcup_G C^g$.
\end{claim}
\begin{proofclaim}
Say $x = ij$ for distinct (possibly commuting) involutions $i, j$. Let $H = C_G^\circ(x)$, a definable, connected, proper subgroup. If $i$ inverts $H$, then let $A = H$; if not, invoke Reineke's Theorem \cite[Corollary~6.5]{BNGroups} and let $A$ be any infinite, definable, connected, subgroup of $C_H^\circ(i) > 1$. In either case $A$ is an infinite, abelian, $\generated{x, i}$-invariant group. We contend that $A$ is contained in a unique conjugate of $C$.

If we fix one Borel subgroup $B$ containing $T$, then by nilpotence $B \leq C$; by conjugacy, all Borel subgroups are contained in conjugates of $C$ (equality is not accessible a priori, since there are no assumptions on the inner structure of $C$). But $A$ certainly extends to a Borel subgroup, hence also to a conjugate of $C$.

Hence we may suppose that $A \leq C$; as the latter is \textsc{ti}, it also is $\generated{x, i}$-invariant. Now there are two cases.
\begin{itemize}
\item
Suppose that $i$ inverts $C$. Let $t \in T \leq C$ and compute:
\[t^x = t^{-xi} = t^{- i x^{-1}} = t^{x^{-1}},\]
so that $x^2$ centralises $t$, and $x^2 \in C(T) = C$.
\item
Suppose not. As there is a $2$-torus $T_i$ containing $i$, with centraliser $C_i$, we see that $C \cap C_i \neq 1$, so by disjunction $C = C_i$. This means that $i \in T$. Bear in mind that $x$ normalises $T$, so that:
\[x^{-1} = x^i = ixi = x i^x i \in x T.\]
This time $x^2 \in T \leq C$.
\end{itemize}
In either case, $x^2$ lies in $C$.
\end{proofclaim}

\begin{claim}
If $x \in G$ is a strongly real element, then $x \in \bigcup_G C^g$.
\end{claim}
\begin{proofclaim}
If $x$ is an involution we are done.
We may suppose that $1 \neq x^2 \in C$, and in particular $x \in N$. If $N/C$ has odd order then we are done. Now suppose that $N/C$ has even order, so Proposition~\ref{p:mainalternative} applies.
Write the smallest definable subgroup containing $x$ as $\defgenerated{x} = \defgenerated{z} \oplus \generated{\alpha}$ with $x = z \alpha$ so that $\defgenerated{z}$ is a $2$-divisible group and $\alpha$ a $2$-element (\cite[Exercise~10 p.~93]{BNGroups} if necessary). Observe how $z$ normalises $C$. Since $N/C$ has order $2$, we find $z \in C$. If $\alpha \in C$ we are done. Otherwise $\alpha \in N \setminus C$ forces $\alpha$ to be an involution inverting $C$, hence also $z$. So $x = z \alpha$ is an involution, which proves $C_x = C_\alpha$.
\end{proofclaim}

Theorem~\ref{t:geometric} now provides the desired contradiction.
This completes the proof of Theorem~\ref{t:hereditary}.
\end{proof}

\subsection{Corollaries of Theorem~\ref{t:hereditary}}

\begin{citedstatement}{Corollary~\ref{c:linear}}
Let $\bK$ be a ranked field of characteristic $0$ and $G \leq \GL_n(\bK)$ be a simple, definable, \emph{non Zariski-closed} subgroup. Then $G$ has no involutions.
\end{citedstatement}
\begin{proof}
Clearly $G$ is $U_2^\perp$. By \cite[Théorème~3]{PQuelques}, Borel subgroups are abelian (hence nilpotent). Moreover, Mustafin \cite[Proposition~2.11]{MStructure} proved conjugacy of the Borel subgroups in all definably linear groups in the ranked setting, so the property is hereditary. Hence we are under the assumptions of Theorem~\ref{t:hereditary}; conclude by simplicity.
\end{proof}

\begin{citedstatement}{Corollary~\ref{c:goodtori}}
Let $G$ be a connected ranked group in which all Borel subgroups are good tori. Then the (possibly trivial) $2$-torsion of $G$ is central.
\end{citedstatement}
\begin{proof}
Here again $G$ is necessarily $U_2^\perp$; moreover, Borel subgroups are maximal good tori, hence abelian and conjugate. And since any definable, connected subgroup of a good torus is still one, the property is hereditary.
\end{proof}

\begin{remarks*}\
\begin{itemize}
\item
It is unclear to us what would happen to Corollary~\ref{c:goodtori} with decent tori instead of good tori since hereditary properties are then lost. The situation is even worse with Cartan subgroups (i.e.~centralisers of maximal decent tori).
\item
Despite vague attempts, we could not prove the following.
\begin{conjecture*}[{Borovik-Burdges \cite[Conjecture~1]{BBDefinably}}]
There is no simple ranked group in which all strongly real elements lie in $\bigcup_G C^g$ for $C$ a Cartan subgroup.
\end{conjecture*}
This does not seem directly related to the $A_1$-conjecture, which we hope to return to shortly.
\end{itemize}
\end{remarks*}

\begin{center}
\rule{.7\textwidth}{.5pt}
\end{center}

This was one of two papers started by the authors in 2018---this is: Paris Album No.~1.  Final details were completed during the trimester program ``Logic and algorithms in group theory''; hospitality of the Hausdorff Institute in Bonn is warmly acknowledged as we are extremely grateful to the organisers and all staff.  We also wish to thank Gregory Cherlin and the anonymous referees for useful comments.

The second author was partially supported on this project by the Sacramento State Research and Creative Activity Faculty Awards Program.

\printbibliography

\end{document}